\documentclass[a4,10.5pt]{article}
\usepackage{amssymb}
\usepackage{amsmath}
\usepackage{latexsym}

\begin{document}

\newcommand{\qed}{\hphantom{.}\hfill $\Box$\medbreak}
\newcommand{\Proof}{\noindent{\bf Proof}\quad}
\newcommand{\vect}[1]{\mbox{\boldmath $#1$}}
\newcommand{\proof}{\noindent{\bf Proof:  }} 

\renewcommand{\theequation}{\thesection.\arabic{equation}} 

\newcommand\one{\hbox{1\kern-2.4pt l }}

\newtheorem{Theorem}{Theorem}[section]
\newtheorem{Lemma}[Theorem]{Lemma}
\newtheorem{Corollary}[Theorem]{Corollary}
\newtheorem{Remark}[Theorem]{Remark}
\newtheorem{Example}[Theorem]{Example}
\newtheorem{Definition}[Theorem]{Definition}

\newtheorem{Construction}[Theorem]{Construction}

\newcounter{Ictr}
\newcommand{\Item}{\refstepcounter{Ictr}\item[(\theIctr)]}

\title{Ideal Secret Sharing Schemes:  
Combinatorial Characterizations, Certain Access Structures, 
and Related Geometric Problems} 

\author{Ryoh Fuji-Hara and Ying Miao\\
Department of Social Systems and Management\\
Graduate School of Systems and Information Engineering\\
University of Tsukuba\\
Tsukuba 305-8573, Japan\\
{\tt \{fujihara, miao\}@sk.tsukuba.ac.jp}}
\date{}
\maketitle

\begin{abstract}
An ideal secret sharing scheme is a method of sharing a secret key 
in some key space among a finite set of participants in such a way 
that only the authorized subsets of participants 
can reconstruct the secret key from their shares 
which are of the same length as that of the secret key. 
The set of all authorized subsets of participants is 
the access structure of the secret sharing scheme. 
In this paper,  we derive several properties and restate 
the combinatorial characterization of an ideal secret sharing scheme 
in Brickell-Stinson model in terms of orthogonality of its representative array. 
We propose two practical models, namely the parallel and hierarchical models, 
for access structures, and then, by the restated characterization, 
we discuss sufficient conditions on finite geometries 
for ideal secret sharing schemes to realize these access structure models.  
Several series of ideal secret sharing schemes realizing 
special parallel or hierarchical access structure model 
are constructed from finite projective planes. 
\end{abstract}

\vskip12pt

{\bf Key words:} Access structure; Combinatorial property; Finite geometry; 
Secret sharing scheme 





\vskip18pt

\begin{center}
{\em Dedicated to Professor Sanpei Kageyama on the Occasion of His Retirement from Hiroshima University}
\end{center}

\vskip18pt


\section{Introduction}
\label{intro}

Secret sharing schemes were first introduced independently by Shamir \cite{Sha79}
and Blakley \cite{BL79} as threshold schemes in 1979.
The general theory of secret sharing schemes was provided 
subsequently by Ito, Saito and Nishizeki \cite{ISN87} in 1987.
Secret sharing schemes become indispensable whenever secret information 
needs to be kept collectively by a group of participants in such a way that 
only a qualified subgroup is able to reconstruct the secret. 
Their applications in information security theory such as 
modeling access control and cryptographic key distribution problems 
were described in \cite{Sim91}. 

Let $\cal{K}$ be a finite key space and $\cal{P}$ be a finite set of participants.
In a secret sharing scheme, a special participant $D \not\in \cal{P}$, called the dealer, 
secretly chooses a key $K \in \cal{K}$ and distributes one share, 
that is, one piece of the shared secret key, 
from the share set $\cal{S}$ to each participant in a secure manner,  
so that no participant knows the shares given to other participants. 

A $(t,n)$ threshold scheme is a secret sharing scheme in which 
if any $t \ (\le n)$ or more participants pool their shares, where $n=|\cal{P}|$,  
then they can reconstruct the secret key $K \in {\cal K}$, 
but any $t-1$ or fewer participants can gain no information about it. 
Threshold schemes have been extensively investigated, 
see, for example, \cite{Miao03, Sti93, SV88}.

In this paper, we consider more general secret sharing schemes.
Let $\Gamma$ be a set of subsets of participants of  ${\cal P}$. 
${\Gamma}$ is called an access structure if 
the subsets in ${\Gamma}$ are exactly those subsets of participants 
that should be able to reconstruct the secret.  
Any subset in ${\Gamma}$ is called an authorized subset of participants.  
Obviously, $\Gamma$ satisfies the monotone property that 
if $T$ is a member of  $\Gamma$, then any subset $R$ of $\cal{P}$ such that
$R \supseteq T$ is also a member of $\Gamma$.

A secret sharing scheme is said to be 
a perfect secret sharing scheme realizing the access structure $\Gamma$, 
if the following two properties are satisfied: 
\begin{itemize}
\item[(1) ]  If the participants of an authorized subset $T \in \Gamma$ 
pool their shares, then they can determine the value of $K \in {\cal K}$.
\item[(2) ]  If the participants of an unauthorized subset $T \not \in \Gamma$ 
pool their shares, then they can determine nothing about the value of $K \in {\cal K}$ 
in an information-theoretic sense, even with infinite computational resources.
\end{itemize}

More precisely, a perfect secret sharing scheme can be represented 
by a set of distributions rules which describe the way 
to distribute secret keys and their corresponding shares.  
A function $f: \{D\} \cup {\cal P} \longrightarrow {\cal K} \cup {\cal S}$ is a distribution rule 
if $f(D) = K \in \cal{K}$ and $f(P) = S \in \cal{S}$  for $P \in \cal{P}$,  
where $f(D)$ is the secret key being shared and 
$f(P)$ is the share given to participant $P \in {\cal P}$. 
Let ${\cal F}$ be the complete set of distributions rules 
from $\{D\} \cup {\cal P}$ to ${\cal K} \cup {\cal S}$, 
and ${\cal F}_K = \{f \in {\cal F}: f(D) = K\}$ for each $K \in {\cal K}$. 
If $K \in {\cal K}$ is the secret key that $D$ wants to share, 
then $D$ will choose uniformly at random a distribution rule $f \in {\cal F}_K$, 
and use it to distribute shares. 

Suppose $\Gamma$ is an access structure and $\cal{F}$ is a set of distribution rules.
It was shown (see \cite{BS92, Sti92}) that $\cal{F}$ is a perfect secret sharing scheme 
realizing the access structure $\Gamma$ if the following two properties are satisfied:
\begin{itemize}
\item[(3)] If there exist $B\in \Gamma$ and  $f,g\in \cal{F}$ such that 
$f(P)=g(P)$ for all $P\in B$, then $f(D)=g(D)$;
\item[(4)] For any unauthorized subset $B \not\in \Gamma$ and  $f\in\cal{F}$,  
there exists a non-negative integer $\lambda(f,B)$  such that for every $K\in\cal{K}$,
$$|\{ g \in {\cal F}_K: g(P) = f(P) \  \textrm{for each}\  P\in B \}| =\lambda(f,B).$$
\end{itemize}
Property (3) means that if a set of participants is authorized, 
then their shares determine the secret key uniquely, 
while property (4) means that if a set of participants is unauthorized, 
then the possibility for any $K \in {\cal K}$ to be the shared secret key is the same.

Although by an obvious modification of the above property (4), 
we can remove the requirement that $f \in {\cal F}_K$ be chosen 
with a uniform probability distribution, in this paper, however, 
we consider only schemes where these probability distributions are uniform. 
 
An important issue in the implementation of secret sharing schemes is the size of shares, 
since the security of such a scheme degrades 
as the amount of the  information that must be kept secret increases. 
Since the secret key $K$ comes from a finite key space ${\cal K}$, 
we can think of $K$ as being represented by a bit-string of length $\log_2 |{\cal K}|$, 
by using a binary encoding, for example.
In a similar fashion, a share can be represented by a bit-string of length $\log_2 |{\cal S}|$.
The information rate of the secret sharing scheme is defined in \cite{BS92} to be 
$\rho = \log_{2}|{\cal K }|/\log_{ 2 }|{\cal S}|$.
It is not difficult to see that the information rate $\rho \le 1$ 
in any perfect secret sharing scheme realizing an access structure $\Gamma$.
In order to be practical, we do not want to have to 
distribute too much secret information as shares. 
Therefore, the information rate should be as close to $1$ as possible.
Since $\rho = 1$ is the optimal situation, 
we refer to such a perfect secret sharing scheme an ideal scheme.

In this paper,  we will investigate ideal secret sharing schemes 
from a combinatorial viewpoint. 
Combinatorial theory provides a very natural setting to model secret sharing schemes.
The three main combinatorial models are 
Brickell-Davenport model, Brickell-Stinson model and entropy model (see \cite{JM98}). 
In the remainder of this paper, 
we will first derive several combinatorial properties and 
restate the combinatorial characterization of an ideal secret sharing scheme in 
Brickell-Stionson model in terms of orthogonality of its representative array 
in Section \ref{prop}. 
We propose two practical models, namely the parallel and hierarchical models, 
for access structures in Section \ref{as}. 
Then by the restated combinatorial characterization, in Section \ref{geop}, 
we discuss sufficient conditions on finite geometries 
for ideal secret sharing schemes to realize these access structure models.  
Several series of ideal secret sharing schemes 
realizing special parallel or hierarchical access structure model 
are constructed from finite projective planes in Section \ref{cons}. 


\section{Combinatorial Properties}
\label{prop}

Let $\Gamma$ be an access structure. 
An authorized subset $B \in \Gamma$ is called minimal 
if $A \not\in \Gamma$ whenever $A \subsetneq B$.
The set of all minimal authorized subsets of $\Gamma$ 
is denoted by $\Gamma^*$, and is called the basis of $\Gamma$.
Since $\Gamma$ consists of all subsets of ${\cal P}$ 
that are supersets of  $B \in \Gamma^*$, 
$\Gamma$ is determined uniquely by $\Gamma^*$, that is,  
$\Gamma = \{C \subseteq {\cal P}:  C \supseteq B, B \in \Gamma^* \}$.
Therefore $B \not\in \Gamma$ if and only if  
$B$ contains no member of $\Gamma^*$ as its subset.
Since $\Gamma^* \subseteq \Gamma$, 
we know that if property (3) is satisfied then the property in which 
 ``$B\in \Gamma$" is replaced by  ``$B \in \Gamma^*$" is also satisfied.
The converse is also true since any $B \in \Gamma$ contains 
at least one member of $\Gamma^*$ as its subset. 

Blundo et al. \cite{BSSV93} depicted a secret sharing scheme 
as an array $M$ (called representative array) of size $|{\cal F}| \times (|{\cal P}| +1)$ 
with entries from ${\cal K} \cup {\cal S}$.
It is supposed that this array is a public knowledge.
The first column of $M$ is indexed by $D$ and consists of elements of ${\cal K}$,
and the remaining columns are indexed by the participants of ${\cal P}$ 
and consist of elements of ${\cal S}$. 
The rows of $M$ correspond to distribution rules of ${\cal F}$. 
$M(f,D) \in {\cal K}$ is the entry in $M$ corresponding to $f \in {\cal F}$ 
and the dealer $D$, and $M(f,P) \in {\cal S}$ is the entry 
corresponding to $f \in {\cal F}$ and $P \in {\cal P}$.
If $B \subseteq {\cal P}$, then $M(f,B)$ denotes $(M(f,P): P \in B)$.
When the dealer $D$ wants to share a secret key $K \in {\cal K}$, 
he will choose uniformly at random a row $f$ of $M$ with $M(f,D) = K$,  
and distribute the share $M(f,P)$ in that row to each participant $P \in {\cal P}$.

With this array representation, 
we can show combinatorial conditions on the array $M$ 
which ensure that the previous conditions (3) and (4) will be satisfied.
\begin{itemize}
\item[(5)] If there exist $B \in \Gamma$ and $f, f' \in {\cal F}$ 
such that $M(f,B) = M(f',B)$, then $M(f,D) = M(f',D)$;
\item[(6)] For any $B \not\in \Gamma$ and $f \in {\cal F}$, 
there exists a non-negative integer $\lambda(f,B)$ such that for any $K \in {\cal K}$,
$$ |\{g \in {\cal F}: M(g,D)=K \textrm{ and } M(g,B)=M(f,B)\}| = {\lambda}(f,B). $$
\end{itemize}
These two conditions will be called ${\cal K}$-uniqueness and ${\cal K}$-balance, 
respectively, in this paper.

Such a combinatorial model was called Brickell-Stinson model 
and the secret sharing scheme was called BS-scheme in \cite{JM98}. 
This model first appeared in \cite{BS92}, 
although a similar idea was first discussed in \cite{BD91}, 
and used in \cite{Mar93, Sti92} and other papers.
If an unauthorized subset of participants pool their shares in a BS-scheme, 
the probability that some $K$ is the secret key is the same 
as that for someone outside the scheme who knows $M$ but none of the shares. 

Since we are interested in ideal secret sharing schemes, and 
a secret sharing scheme is ideal if and only if $|{\cal K}| = |{\cal S}|$, 
from now on we assume  ${\cal K} = {\cal S}$.

Based on the basic results described above, in the rest of this section, 
we try to derive a few combinatorial properties of 
ideal perfect secret sharing schemes in Brickell-Stinson model 
in terms of the orthogonality of their corresponding arrays $M$. 

Le $A$ be an $N \times m$ array with entries from a set $S$ of $s$ symbols.
Consider an $N \times t$ sub-array $A_0$ consisting of $t$ columns of $A$.
If every $t$-tuple of $S^t$ appears exactly $\lambda$ times in $A_0$ as a row, 
then these $t$ columns are said to be orthogonal.  
Obviously $\lambda = N / s^t$.
If any $t$ columns of $A$ are orthogonal, 
then this array $A$ is an orthogonal array of strength $t$,  
denoted by $OA_{\lambda}(N,m,s,t)$. 
Clearly, if $A$ is an orthogonal array of strength $t$, 
then it is also of strength $t'$ for $1 \le t' \le t$.
Let $A'$ be an $N \times n$ sub-array consisting of certain $n$ columns of $A$.  
If the maximum number of orthogonal columns in $A'$ is $k$, 
then $k$ is called the orthogonal rank of $A'$, denoted by $\phi(A')$.
Under the assumption that $\phi(A') = k < n$, if for any $k$ orthogonal columns 
$c_1, c_2, \ldots, c_k$ of $A'$, and for any two distinct rows 
$(a_1,a_2, \ldots, a_n)$, $(b_1,b_2,\ldots, b_n)$ of $A'$,  
the condition $a_{c_1} = b_{c_1}, \ldots, a_{c_k} = b_{c_k}$ can imply 
$a_{c_{k+1}}=b_{c_{k+1}}, \ldots, a_{c_n} = b_{c_n}$, 
then $A'$ is said to be regular. 

In an ideal secret sharing scheme, 
if $\Gamma^*$ consists of all $k$-subsets of ${\cal P}$, 
then the scheme is an ideal  $(k, |{\cal P}|)$ threshold scheme.
Ideal threshold schemes were characterized by Jackson and Martin \cite{JM98} 
in terms of orthogonal arrays.

\begin{Theorem}
There exists an ideal $(k,n)$ threshold scheme with $|{\cal K}| = s$ if and only if 
there exists an orthogonal array $OA_1(s^k, n+1, s, k)$.
\end{Theorem}

However, the relationship between an ideal secret sharing scheme 
realizing some arbitrary access structure and 
the orthogonality of the $|{\cal F}| \times (|{\cal P}| +1)$ representative array $M$ 
used to depict it has not been made clear yet.
We try to ascertain the relationship between them.
We always assume that every participant belongs to some authorized subsets, 
and not every $B \in {\Gamma^*}$ consists of only one participant,  
otherwise the scheme would be degenerate, in which we are not interested.                                       
Let $M(B)$ be the $|{\cal F}| \times |B|$ sub-array of $M$ 
consisting of columns of $B \subseteq {\cal P}$. 

\begin{Lemma}
\label{freq}
In the column indexed by $D$ and in any column of $M$ indexed by 
a participant belonging to some subset of ${\Gamma^*}$, 
every element of ${\cal S} = {\cal K}$ appears equally often.
\end{Lemma}

\proof
Without loss of generality, we may assign an arbitrary order to each participant of ${\cal P}$. 
For any ordered subset $A \subseteq {\cal P}$, 
let $R_A(\vect{x})$ denote the set of rows of $M$ having 
the $|A|$-tuple $\vect{x} \in {\cal S}^{|A|}$ 
in the columns indexed by the participants in $A$. 

We choose a subset $B \in {\Gamma^*}$ with $|B| \ge 2$. 
Let $C = B \setminus \{P\}$, where $P \in B$.  
Then $C \not\in {\Gamma^*}$, and hence from the $\cal{K}$-balance condition, 
each element of ${\cal K} = {\cal S}$ appears exactly $|R_C(\vect{x})|/|{\cal K}|$ times 
in the multiset $\{M(r,D): r \in R_C(\vect{x})\}$ for any $\vect{x} \in {\cal S}^{|C|}$ 
in the columns indexed by $C$.  
Therefore every element of ${\cal K} = {\cal S}$ appears 
${\sum}|R_C(\vect{x})|/|{\cal K}| = |{\cal F}|/|{\cal K}|$ times in the column indexed by $D$, 
where the summation is over the set of all $|C|$-tuples $\vect{x} \in {\cal S}^{|C|}$ 
occurred in the columns indexed by $C$. 

Now we consider the column indexed by $P \in B$.
For any $\vect{x} \in {\cal S}^{|C|}$ in the columns of $C$, 
since $B \in  {\Gamma^*}$ satisfies the ${\cal K}$-uniqueness condition, 
there exists a one-to-one mapping ${\xi}_{\vect{x}}$ from ${\cal K}$ to ${\cal S}$
which maps an element $K \in {\cal K}$ in column $D$ 
to an element $S \in {\cal S}$ in column $P$  
in the same row of $R_C(\vect{x})$. 
We have already known that each element of ${\cal K}$ appears 
exactly $|R_C(\vect{x})|/|{\cal K}|$ times in the multiset $\{M(r,D): r \in R_C(\vect{x})\}$, 
so we can also know that each element of ${\cal S}$ appears 
exactly $|R_C(\vect{x})|/|{\cal S}|$ times in the multiset $\{M(r,P): r \in R_C(\vect{x})\}$. 
Altogether, every element of ${\cal S}$ appears 
${\sum}|R_C(\vect{x})|/|{\cal S}| = |{\cal F}|/|{\cal S}|$ times in the column indexed by $P$, 
where the summation is again over the set of all $|C|$-tuples $\vect{x} \in {\cal S}^{|C|}$ 
occurred in the columns indexed by $C$. 

Also, from the ${\cal K}$-uniqueness condition, we can easily know that 
every element of ${\cal S} = {\cal K}$ appears equally often 
in any column of $M$ indexed by a participant belonging to 
some singleton subset of ${\Gamma^*}$.
\qed

\begin{Lemma}
\label{baseortho}
If $B \in \Gamma^*$, then $\phi(M(B)) = |B|$.
\end{Lemma}

\proof
We prove this lemma by induction. 
From Lemma \ref{freq}, each column of $B$ is orthogonal.  
Let $C$ be an ordered $t$-set of columns in $B$ 
which are orthogonal where $t = |C| \le |B|-1$.
For any $\vect{x} \in {\cal S}^t$, 
$R_C(\vect{x})$ consists of $|{\cal F}|/|{\cal S}|^t$ rows.  
Let $R_C(K,\vect{x})$ denote the set of rows in $R_C(\vect{x})$ 
having $K \in {\cal K} = {\cal S}$ in column $D$. 
Since $C \not\in {\Gamma}^*$ satisfies the ${\cal K}$-balance condition,
$|R_C(K,\vect{x})| = |{\cal F}|/|{\cal S}|^{t+1}$ for any $K \in {\cal K} = {\cal S}$.
Let $\vect{a}$ be a column in $B \setminus C$ and let $E = B \setminus (C \cup \{\vect{a}\})$. 
If $t \le |B|-2$, then $E \neq \emptyset$. 
Suppose that $\vect{y} \in {\cal S}^{|E|}$ appears 
in a row of $R_{C}(\vect{x})$ in the columns of $E$. 
Consider the set of rows $R_C(\vect{x}) \cap R_E(\vect{y})$ in $M$, 
denoted by $R_{C \cup E}(\vect{x}||\vect{y})$, where $\vect{x}||\vect{y}$ 
means the concatenation of $\vect{x}$ and $\vect{y}$ with the assigned order.  
Since $C \cup E \not\in {\Gamma}^*$ satisfies the ${\cal K}$-balance condition, 
every element $K$ of ${\cal K} = {\cal S}$ occurs as a symbol in column $D$ 
in the rows of $R_{C \cup E}(\vect{x}||\vect{y})$ 
exactly $|R_{C \cup E}(\vect{x}||\vect{y})|/|{\cal S}|$ times
and in the rows $R_C(\vect{x})$ exactly 
$|R_C(\vect{x})|/|{\cal S} = |{\cal F}|/|{\cal S}|^{t+1}$ times. 
Since $B = C \cup E \cup \{\vect{a}\} \in {\Gamma}^*$ satisfies 
the ${\cal K}$-uniqueness condition, there is 
a one-to-one mapping $\xi_{\vect{x}||\vect{y}}$,  in $R_{C \cup E}(\vect{x}||\vect{y})$, 
from ${\cal K}$ of column $D$ to ${\cal S}$ of column $\vect{a}$.
Therefore, the total number of symbol $S \in {\cal S}$ occurs in column $\vect{a}$ 
in the rows of $R_C(\vect{x})$ is $|R_C(\vect{x})/|{\cal S}| = |{\cal F}|/{\cal S}^{t+1}$, 
which is the same as that in column $D$.
This implies that the columns of $C \cup \{\vect{a}\}$ are orthogonal.
If $t = |B|-1$, then $E = \emptyset$, that is, $B = C \cup \{\vect{a}\}$, 
and there is a one-to-one mapping $\xi_{\vect{x}}$,  in $R_{C}(\vect{x})$, 
from ${\cal K}$ of column $D$ to ${\cal S}$ of column $\vect{a}$.
In a similar fashion, we can show that the columns of $C \cup \{\vect{a}\}$ are orthogonal.
Therefore, the columns of $B$ are orthogonal.
\qed

The following is a combinatorial characterization of an ideal secret sharing scheme.

\begin{Theorem}
\label{char}
An array $M'$ represents an ideal secret sharing scheme 
realizing an access structure $\Gamma$ 
if and only if $M'$ satisfies the following conditions:
\begin{itemize}
\item[(7)] For any $B \in \Gamma^*$,  $M'(\{D\}\cup B)$ is a regular orthogonal array 
of strength $|B|$ but not $|B|+1$;
\item[(8)] For any $B \not\in \Gamma$, $\phi(M'(\{D\} \cup B)) = \phi(M'(B)) +1$ 
and $M'(\{D\} \cup B)$ is regular.
\end{itemize}
\end{Theorem}

\proof We first consider the necessity. Suppose $M'$ represents 
an ideal secret sharing scheme realizing an access structure $\Gamma$.  
If $B \in {\Gamma}^*$, then by Lemma \ref{baseortho}, 
we know that the columns of $B \in {\Gamma^*}$ are orthogonal.
Also, for any column $\vect{a}$ in $B \in {\Gamma}^*$,  by the proof of Lemma \ref{freq}, 
the columns of $\{D\} \cup (B \setminus \{\vect{a}\})$ are orthogonal. 
Furthermore, since $B \in {\Gamma^*}$, 
there is a one-to-one mapping from ${\cal K} = {\cal S}$ of column $D$ 
onto ${\cal S}$ of column $\vect{a}$.
Therefore $M'(\{D\}\cup B)$ is a regular orthogonal array of strength $|B|$ but not $|B|+1$. 
If $B \not\in {\Gamma}$, then for any $C \subseteq B$, we know that $C \not\in {\Gamma}^*$. 
We define $O \subseteq B$ to be the set of $\phi(M'(B))$ orthogonal columns. 
Since $O \not\in {\Gamma}^*$, from the perfectness of the scheme, 
we immediately know that the columns of $\{D\} \cup O$ should be also orthogonal, 
and furthermore, $\phi(M'(\{D\} \cup B)) = \phi(M'(B)) +1$.  

Conversely, if there exists $B \subseteq {\cal P}$ 
such that $B \in \Gamma$ and $M'(r,B) = M'(r',B)$, 
then there exists $B' \subseteq B$ such that $B' \in {\Gamma}^*$ and $M'(r,B') = M'(r',B')$. 
By (7), we obtain that $M'(r,D) = M'(r',D)$. 
For any $B \not\in \Gamma$, 
since  $\phi(M'(\{D\} \cup B)) = \phi(M'(B)) +1$, 
if we assume $O \subseteq B$ to be the set of $\phi(M'(B))$ orthogonal columns, 
then columns in $\{D\} \cup O$ are orthogonal.  
Therefore for any distinct $K, K' \in {\cal K}$ and for any $\vect{x} \in {\cal S}^{\phi(M'(B))}$,
$$ |\{r: M'(r,D)=K \textrm{ and } M'(r,O)=\vect{x}\}| = |\{r: M'(r,D)=K' \textrm{ and } M'(r,O)=\vect{x}\}|. $$
Since $M'(\{D\} \cup B)$ is regular, 
we know that for any distinct $K, K' \in {\cal K}$ and for any $\vect{y} \in {\cal S}^{|B|}$,
$$ |\{r: M'(r,D)=K \textrm{ and } M'(r,B)=\vect{y}\}| = |\{r: M'(r,D)=K' \textrm{ and } M'(r,B)=\vect{y}\}|. $$
\qed

Note that similar but a little different results can be found in \cite{BD91}, 
and characterizations of ideal secret sharing schemes in different terms 
can be found in, for example, \cite{BS92, JM98, MFP03}.    

In the following sections, we will investigate 
ideal secret sharing schemes with certain specified access structures 
by means of finite geometries via Theorem \ref{char}.    
Therefore, now we turn to consider a vector space of dimension $t$ over a finite field $GF(q)$, 
denoted by $GF(q)^k$.
Let $G$ be the $k \times m$ matrix defined by 
$$ G = (\vect{v}_1, \vect{v}_2, \ldots, \vect{v}_m ), \ \ \ \ \ 
\vect{v}_i \in GF(q)^k \setminus \{(0,0,\ldots,0)^T\}, $$
and $M''$ be the $q^k \times m$ array generated by $G$
$$M'' = (\vect{x}^T G :  \vect{x} \in GF(q)^k),$$
where $\vect{x}^T$ denotes the transpose of  $\vect{x} \in GF(q)^k$. 
Then we have the following well-known result (see, for example, \cite{HSS99}):

\begin{Theorem}
\label{linorth}
The vectors $\vect{v}_{i_1}, \vect{v}_{i_2}, \ldots, \vect{v}_{i_t}$, where $t \le k$, 
are linearly independent if and only if 
the $i_1$th, $i_2$th, $\ldots$, $i_t$th columns of $M''$ are orthogonal.
\end{Theorem}

By Theorem \ref{linorth}, we know that if any $t$ vectors of $G$ are linearly independent 
then $M''$ is an orthogonal array $OA_{\lambda}(q^k, m, q, t)$ with $\lambda = q^{k-t}$.
Since the set of rows of $M''$ is a vector subspace of $GF(q)^m$, 
such an orthogonal array is usually called a linear orthogonal array.

\begin{Corollary}
\label{M''}
$M''$ is a regular linear orthogonal array of strength at least one.
\end{Corollary}

\Proof
Since $M''$ is generated by $G = (\vect{v}_1, \vect{v}_2, \ldots, \vect{v}_m )$, and   
$\vect{v}_i \neq (0,0,\ldots,0)^T$ for $1 \le i \le m$, 
any vector $\vect{v}_i$ is linearly independent, 
and thus $M''$ is a linear orthogonal array of strength at least one. 
Now we consider its regularity. 
Let $X$ be a $k \times n$ sub-matrix of $G$ consisting of certain $n$ columns of $G$, 
say $X = (\vect{x}_1, \vect{x}_2, \ldots, \vect{x}_n )$, with rank $t \ ( < n)$. 
Assuming $\vect{x}_{i_1}, \vect{x}_{i_2}, \ldots, \vect{x}_{i_t}$ are linearly independent. 
Then each of the remaining vectors $\vect{x}_{i_{t+1}}, \vect{x}_{i_{t+2}}, \ldots, \vect{x}_{i_n}$ 
can be represented by a linear combination of $\vect{x}_{i_1}, \vect{x}_{i_2}, \ldots, \vect{x}_{i_t}$. 
This implies that for any sub-row $(a_1, a_2, \ldots, a_n)$ of $M''$ indexed by $X$, 
$a_{i_{t+1}}, a_{i_{t+2}}, \ldots, a_{i_n}$ can be uniquely determined by 
$a_{i_1}, a_{i_2}, \ldots, a_{i_t}$, 
which shows the regularity of $M''$.
\qed

In the description of the relationship between ideal secret sharing schemes 
and orthogonal arrays, regularity is an important character.
Kuriki \cite{Kur04} showed a construction for non-regular orthogonal arrays.

If we assume that the dealer and participants are assigened vectors of $GF(q)^k$, that is 
$D = \vect{d}$, ${\cal P} = \{\vect{v}_1, \vect{v}_2, \ldots, \vect{v}_m\}$, 
where $\vect{d}, \vect{v}_i \in GF(q)^k$,  
then the following theorem is an immediate consequence of 
Theorem \ref{char}, Theorem \ref{linorth},  and Corollary \ref{M''}. 

\begin{Theorem}
\label{vectchar}
An array $M''$ generated by $G = (\vect{d}, \vect{v}_1, \vect{v}_2, \ldots, \vect{v}_m )$
represents an ideal secret sharing scheme 
realizing an access structure $\Gamma$ if and only if 
$M''$ satisfies the following conditions: 
\begin{itemize}
\item[(9)]  For any $X \in \Gamma^*$, $|X| \le k$,  any $|X|$ vectors of $\{\vect{d}\} \cup X$ 
are linearly independent but the $|X|+1$ vectors of $\{\vect{d}\} \cup X$  are linearly dependent;
\item[(10)]  For any $Y \not\in \Gamma$,   $rank(\{\vect{d}\} \cup Y) = rank(Y) + 1 \le k$.
\end{itemize}
\end{Theorem}

The reader is referred to \cite{MvD97} for the characterizations of linear secret sharing schemes 
to compare Theorem~\ref{vectchar} with the related part in \cite{MvD97}.


\section{Access Structures}
\label{as}

In Section \ref{intro}, we defined the access structure $\Gamma$ 
of a secret sharing scheme as the set of all authorized subsets 
of the participants ${\cal P}$ which satisfies the monotone property that 
if $T$ is a member of $\Gamma$, then any subset of ${\cal P}$ containing $T$ 
as a subset is also a member of  $\Gamma$.
We also denoted the set of all minimal authorized subsets of $\Gamma$, 
or the basis of $\Gamma$, by $\Gamma^*$.
Since $\Gamma$ is determined uniquely by $\Gamma^*$, in order to 
define $\Gamma$, it is enough to define its minimal authorized subsets $\Gamma^*$. 

However, in some specified access structure models, 
not all given sets of subsets of ${\cal P}$ can be used to define $\Gamma^*$,
even if these subsets are mutually non-included.   
In order to meet some specified requirements of an access structure model, 
it may happen that some subsets are not allowed to be in $\Gamma^*$ 
but some others have to be added to $\Gamma^*$.
It is understood that it is not easy to define suitable $\Gamma^*$ 
which can be realized by some ideal secret sharing scheme.

Uehara et. al. \cite{UNN86} and Jackson and Martin \cite{JM98} 
showed some relationship between ideal secret sharing schemes and matroids.
We do not discuss about metroids in the paper. 
Instead, we investigate the following two access structure models  
which can be realized by ideal secret sharing schemes.

\subsection{Parallel Model} 

Suppose that the set ${\cal P}$ of participants is partitioned into 
mutually disjoint groups $G_1, G_2, \ldots, G_m$. 
These groups can be considered as the departments of a university or a company.
We may want to set up a committee to manage such affairs as 
Ph.D.  degree defenses or personnel affairs, 
and we expect the members of the committee to meet certain requirements such as: 
(i) the members should consist of at least $t$ persons; 
(ii) the members should not be all from one department. 

We consider the following model. 
A $t$-subset of the participants ${\cal P}$ is defined to be in the basis $\Gamma^*$ 
of the access structure if and only if the $t$-subset is not included 
in any group of $G_1, G_2, \ldots, G_m$, that is,
$$ \Gamma^* = {{\cal P} \choose t} \setminus \bigcup\limits_{i=1}^m { G_i \choose t}, $$
where, for a set $S$ and an integer $t$, ${S \choose t}$ means 
the set of all $t$-subsets of $S$.
Such a model is called a parallel model of strength $t$ in this paper.

\subsection{Hierarchical Model}

In this subsection, we consider another model. 
We suppose that the set ${\cal P}$ of participants is composed of two disjoint groups,  
the upper group ${\cal U}$ and the lower group ${\cal L}$. 
We may think of ${\cal U}$ as a group of managers and 
${\cal L}$ a group of lower employees in some organization.  
A contract should be signed by at least $t$ participants in this organization where 
at least one of the signers should come from the upper group ${\cal U}$.

A $t$-subset of the participants ${\cal P}$ is defined to be in the basis $\Gamma^*$ 
of the access structure if and only if 
at least one member of the $t$-subset is from the upper group, that is,
$$ \Gamma^* = {{{\cal U} \cup {\cal L}} \choose t} \setminus {{\cal L} \choose t}. $$
We call this model a hierarchical model of strength $t$.

Similar models and the reliability of the schemes were 
considered and analyzed in \cite{Koy81, Bric89}. 
Researches on secret sharing schemes on some special access structures 
can also be found in, for example, \cite{PS00, MFP04, DS06, Mar07}. 
A survey of geometrical contributions to secret sharing theory is provided in \cite{JMO04}. 
In the remainder of this paper, we focus our attention on the construction of 
ideal secret sharing schemes which can realize the two access structure models 
defined above by using finite geometries.


\section{Geometric Problems}
\label{geop}

Let $\Pi_r(q)$ (or $\Pi_r$ if no confusion arises) be a finite projective space 
of dimension $r$ over $GF(q)$. 
It is well known that $\Pi_{r}(q)$ can be constructed from all non-zero $(r+1)$-tuples 
of elements of $GF(q)$ by regarding the points as the equivalence classes $[\vect{x}]$ 
where the equivalence relation $\vect{x} \sim \vect{y}$ is defined to be 
$\vect{x} = {\lambda}\vect{y}$ for some non-zero ${\lambda} \in GF(q)$.  
It can be easily seen that a set $S$ of $t$ points in $\Pi_r$, $r \ge  t-1$, 
is linearly independentif and only if there is no subspace $\Pi_{t-2}$ which contains $S$.  
Clearly, any two points in a projective space $\Pi_r$, $r \ge 1$, is always independent.
Let ${\cal K}$ be a set of $k \ (\ge t)$ points in $\Pi_r$.
If any set of $t$ points of ${\cal K}$ is independent, then we call ${\cal K}$ a $t$-independent set.
A $k$-arc in $\Pi_r$ is a set ${\cal K}$ of $k$ points which is an $(r+1)$-independent set.
To find the maximum size of $t$-independent sets in $\Pi_r$ is a very difficult problem.
The interesting reader is referred to Hirschfeld \cite{Hir98} for several known results.

A pencil in a projective plane $\Pi_2$ is the set of lines through some common point 
which lie in this plane.  
In general, a pencil is the set of hyperplanes $\Pi_{r-1}$ 
containing a given $\Pi_{r-2}$ in common and lying in $\Pi_{r}$.

Let $\Psi_0, \Psi_1, \ldots, \Psi_q$ be a pencil in $\Pi_{r}(q)$. 
A pencil arc ${\cal K}$ is a set of $k$ points, denoted by $k$-parc,  
in $\Pi_{r}(q)$ satisfying the following conditions:
\begin{enumerate}
\item[(11)] Each ${\cal K} \cap \Psi_i$ is an $h_i$-arc in $\Pi_{r-1}$ for $0 \le i \le q$, 
where $h_i = |{\cal K}  \cap  \Psi_i|$;
\item[(12)] ${\cal K} \cap \Psi_i \cap \Psi_j = \emptyset$ for $i \ne j$, where $0 \le i, j \le q$; 
\item[(13)] Any $r+1$ points of ${\cal K}$ not contained in any single $\Psi_i$ are independent.
\end{enumerate}
A regular pencil arc, denoted by  $(m,h)$-parc, is a pencil arc in which 
$h_0 = h_1 = \ldots = h_{m-1} = h$ and $h_m = h_{m+1} = \ldots = h_q = 0$. 

\begin{Theorem}
Let $h_i \ge 1$ for $0 \le i \le m-1$ and $h_0 = \min\{h_i: 0 \le i \le m-1\}$. 
If there exists a $k$-parc in $\Pi_{t-1}$ with $k = h_0 + h_1 + \ldots + h_{m-1}$ points, 
then there exists an ideal secret sharing scheme
of parallel model with $k - h_0$ participants and of strength $t$.
\end{Theorem}

\Proof
Let ${\cal K}$ be a $k$-parc in ${\Pi}_{t-1}(q)$ 
with $\Psi_0, \Psi_1, \ldots, \Psi_q$ being its pencil, 
and ${\cal K}_i = {\cal K} \cap \Psi_i$ with $h_i = |{\cal K}_i|$ for $0 \le i \le q$.
We assign the dealer $D$ an arbitrary point in ${\cal K}_0$ and 
each of the participants in ${\cal P}$ a different point in 
${\cal K}_1 \cup {\cal K}_2 \cup \ldots \cup {\cal K}_{m-1}$. From the definition of a $k$-parc, 
every set of $t-1$ points in ${\cal K}_i$ is independent, 
but no set of $t$ points in ${\cal K}_i$ is independent for $0 \le i \le m-1$, 
nor set of $t+1$ points in $\Pi_{t-1}$ is independent, 
and every set  $T$ of $t$ points not contained in any ${\cal K}_i$, i.e.,  
$T \not\subseteq {\cal K}_i$ for $0 \le i \le m-1$, is independent.
Therefore, any set of $t+1$ points consisted of $D$ and the other 
$t$ points from ${\cal K}_1 \cup {\cal K}_2 \cup \ldots \cup {\cal K}_{m-1}$ 
which are not all included in any ${\cal K}_i$, $0 \le i \le m-1$, is a $t$-independent set. 
Also, any set of $t+1$ points consisted of $D$ and the other 
$t$ points all from some ${\cal K}_i$, $1 \le i \le m-1$, is not independent, 
but $rank(\{D\} \cup Y) = rank(Y) + 1 = t$ 
for any set $Y$ of $t$ points completely contained in some ${\cal K}_i$, $1 \le i \le m-1$.    
Therefore, the two conditions of Theorem \ref{vectchar} are satisfied. 
\qed

Now we define another type of geometric structure. 
Let $\Psi$ be a hyperplane of $\Pi_{r}$, 
${\cal K}_{1}$ be a set of $k_1$ points in $\Pi_{r} \setminus \Psi$,  
and ${\cal K}_2$ be a set of $k_2$ points in $\Psi$.
A hierarchical arc in $\Pi_{r}$ is a set ${\cal K} =  {\cal K}_1 \cup {\cal K}_2$ 
of $k_1 + k_2$ points in $\Pi_{r}$, denoted by $(k_1, k_2)$-harc, 
satisfying the following conditions:
\begin{enumerate}
\item[(14)] ${\cal K}_1$ is a $k_1$-arc in $\Pi_{r}$;  
\item[(15)] ${\cal K}_2$  is a $k_2$-arc in $\Pi_{r-1}$; 
\item[(16)] Any $r+1$ points of ${\cal K}$ not contained in the hyperplane $\Psi$ 
are independent.
\end{enumerate}

\begin{Theorem}
Let $k_1 \ge 2$ and $k_2 \ge 0$. 
If there exists a $(k_1, k_2)$-harc in $\Pi_{t-1}$,  
then there exists an ideal secret sharing scheme of hierarchical model 
with $k_1+k_2-1$ participants and of strength $t$.
\end{Theorem}

\Proof
Let ${\cal K} = {\cal K}_1 \cup {\cal K}_2$ be a $(k_1, k_2)$-harc in ${\Pi}_{t-1}$ 
with ${\cal K}_1$ being a $k_1$-arc in $\Pi_{t-1}$ 
and ${\cal K}_2$ being a $k_2$-arc in $\Pi_{t-2}$.
We assign the dealer $D$ an arbitrary point in ${\cal K}_1$ and 
each of the participants in ${\cal P}$ a different point in ${\cal K} \setminus \{D\}$. 
Let ${\cal K}_1 \setminus \{D\}$ be the upper group, and ${\cal K}_2$ be the lower group.
From the definition of a $(k_1, k_2)$-harc, 
every set of $t+1$ points consisted of $D$ and 
the other $t$ points from a set $X \subseteq {\cal K} \setminus \{D\}$,
where $X$ contains at least one point of ${\cal K}_1 \setminus \{D\}$, 
is not independent, 
but any of its $t$-subsets of points is independent. 
Also, any set of $t+1$ points consisted of $D$ and the other 
$t$ points all from ${\cal K}_2$ is not independent, 
but $rank(\{D\} \cup Y) = rank(Y) + 1 = t$ 
for any set $Y$ of $t$ points completely contained in ${\cal K}_2$.    
Therefore, the two conditions of Theorem \ref{vectchar} are satisfied. 
\qed

In the remainder of this section, we consider the upper bounds 
for existence of a $k$-parc and a $(k_1,k_2)$-harc in $\Pi_2(q)$.

\begin{Theorem}
\label{parcb2}
If there exists a $k$-parc in $\Pi_2(q)$ then $k \le 2q$.
\end{Theorem}

\Proof
Let ${\cal K}$ be a $k$-parc in ${\Pi}_{2}(q)$ with $\Psi_0, \Psi_1, \ldots, \Psi_q$ being its pencil,  
${\cal K}_i = {\cal K} \cap \Psi_i$ with $h_i = |{\cal K}_i|$ for $0 \le i \le q$, 
and $h_0 = \max\{h_i: 0 \le i \le q\}$.
Let $l$ be a line through a point in ${\cal K}_0$ but $l \neq \Psi_0$.
Then $l$ can meet ${\cal K} \setminus  \Psi_0$ in at most one point. 
Therefore, $\sum_{i=1}^{n}{h_i} \le q$, 
which implies that $k = \sum_{i=0}^{n}{h_i} \le q = h_0 + \sum_{i=1}^{n}{h_i} \le 2q$.
\qed

The bound in Theorem \ref{parcb2} is attained 
when ${\cal K}$ has $q$ points on two lines $\Psi_0$ and $\Psi_i$, $1 \le i \le q$, 
respectively, but no points on any other lines. 
We do not know any other examples of $k$-parcs with the maximal number of points.

\begin{Theorem}
\label{harcb2}
If there exists a $(k_1, k_2)$-harc in $\Pi_2(q)$ then $k_1 + k_2 \le q+2$.
\end{Theorem}

\Proof
Let ${\cal K} = {\cal K}_1 \cup {\cal K}_2$ be a $(k_1, k_2)$-harc in ${\Pi}_{2}(q)$ 
with $\Psi$ being its underlying hyperplane, 
${\cal K}_1$ being a $k_1$-arc in $\Pi_2$ and ${\cal K}_2$ being a $k_2$-arc in $\Pi_{1}$. 
Let $P \in {\cal K}_1$ and $l$ be a line through the point $P$. 
Then $l$ can meet ${\cal K} \setminus  \{P\}$ in at most one point. 
Therefore, $(k_1-1)+k_2 \le q+1$, 
which implies that $k_1+k_2 \le q + 2$.
\qed

A $k$-parc (or a $(k_1, k_2)$-harc) is said to be maximal 
if the number of points it contained 
attains the above bound in Theorem \ref{parcb2} (or Theorem \ref{harcb2}, respectively). 


\section{Constructions of Pencil or Hierarchical Arcs}
\label{cons}

In this section, we focus on the constructions of 
ideal secret sharing schemes of parallel or hierarchical model with strength $3$. 
To this end, we study $k$-parcs and $(k_1, k_2)$-harcs 
in a finite desarguesian projective plane $\Pi_2(q)$.

Obviously, if there exists a $k$-arc ${\cal K}$ in $\Pi_2(q)$, 
then there exists a $k$-parc ${\cal K}$ with pencil $\Psi_0, \Psi_1, \ldots, \Psi_q$ such that 
$|{\cal K} \cap {\Psi_i}| =  0,  1 \ \mbox{or} \ 2$ for each line $\Psi_i$, $0 \le i \le q$, 
and also there exists a $k$-harc 
${\cal K} = ({\cal K} \cap {\Psi}) \cup ({\cal K} \setminus ({\cal K} \cap {\Psi}))$, 
where ${\Psi}$ is a hyperplane in ${\Pi}_2(q)$,  
such that $|{\cal K} \cap {\Psi}| = 1 \ \mbox{or} \ 2$
(note that $|{\cal K} \cap {\Psi}| = 2$ when $q$ is even). 
In this section we are more interested in pencil arcs and hierarchical arcs 
with more points in $\Psi_i$ and $\Psi$, respectively.

Let us recap some definitions and results in finite projective geometry. 
For terms undefined and results mentioned without references,  
the reader is referred to \cite{FV83, FJ84, Hir98} for the detailed information. 
Let $\Pi$ be a projective plane with point set ${\cal P}$ and line set ${\cal L}$.
A subplane ${\Pi}'$ of $\Pi$ is a projective plane with point set ${\cal P}'$ and ${\cal L}'$ 
where ${\cal P}' \subseteq {\cal P}$ and ${\cal L}' \subseteq {\cal L}$.
If the order $q'$ of ${\Pi}'$ is less than the order $q$ of $\Pi$, then $q' \le \sqrt{q}$.  
A subplane $\Pi'$  of order $\sqrt{q}$ in a projective plane $\Pi$ of order $q$ 
is called a Baer subplane of $\Pi$. 
A Baer subplane ${\Pi}'$ intersects each line of ${\Pi}$ in 
either $\sqrt{q}+1$ points or one point.  
Any point outside ${\Pi}'$ is on exactly one line of ${\Pi}'$. 
A projective plane $\Pi$ has a Baer subplane if and only if the order of $\Pi$ is a square.

Let $\Pi$ be a projective plane of order $q^2$ and $l_\infty$ be an arbitrary line of $\Pi$.
Consider a Baer subplane ${\Pi}'$ of $\Pi$ which intersects $l_\infty$ in $q+1$ points, 
that is, ${\Pi}' \cap l_{\infty}$ is a line of ${\Pi}'$. 
For any line $l \neq l_{\infty}$ of ${\Pi}$ such that 
${\Pi}' \cap l$ is a line of ${\Pi}'$ consisting of $q+1$ points, 
the line ${\Pi}' \cap l$ of ${\Pi}'$ is called a segment of the line $l$ of ${\Pi}$,
whereas $({\Pi}' \cap l) \setminus l_{\infty}$ 
consisting of $q$ points in the affine part $\Pi \setminus l_\infty$ 
is called an affine segment.
$\Pi' \setminus l_\infty$ is called an affine Baer subplane 
of the affine plane $\Pi \setminus l_\infty$.

In the remainder of this section, 
we always suppose that $\Pi$ is a projective plane of order $q^2$  
and ${\Pi}'$ a Baer subplane of $\Pi$. 
Let $l_1, l_2, \ldots, l_{q}$ and $l_\infty$ be $q+1$ lines of $\Pi$ such that 
$\{{\Pi}' \cap l_1, \ldots, {\Pi}' \cap l_{q}, {\Pi}' \cap l_\infty\}$ is a pencil of $\Pi'$ 
with a common point $P \in {\Pi}'$.
Let $\cal F$ be the set of Baer subplanes of $\Pi$ 
each of which intersects each of $l_1, l_2, \ldots, l_{q}$ and $l_\infty$ in $q+1$ points.
The following known results on Baer subplanes can be found, for example, 
in \cite{ FV83, FJ84, Hir98}. 

\begin{enumerate}
\item[(17)]  Each affine Baer subplane of $\cal F$ intersects each $l_i \setminus \{P\}$, 
$i = 1, \ldots, q, \infty$,  in an affine segment.  
For each $i$, $i = 1, \ldots, q, \infty$, the set of affine segments is denoted by ${\cal L}_i$. 
\item[(18)] For each $i$, $i = 1, \ldots, q, \infty$, the  incidence structure 
$\Theta_i = ( l_i \setminus \{P\}, {\cal L}_i)$ is isomorphic to an affine plane of order $q$.
\item[(19)]  For any line $l$ of $\Pi$ not containing $P$,  
the intersections of $l$ with $l_1, l_2, \ldots, l_q$ form an affine segment. 
Let ${\cal L}$ denote the set of such affine segments.
\item[(20)]  The incidence structure of the points in 
$l_1\cup l_2 \cup \ldots \cup l_q \setminus \{P\}$ and the affine segments 
in ${\cal L}_1 \cup {\cal L}_2 \cup \ldots \cup {\cal L}_q \cup {\cal L}$ 
forms an affine three space of order $q$, denoted by $\Lambda$. 
\item[(21)] There are $q$ mutually disjoint affine Baer subplanes 
which lie on $l_1, l_2, \ldots, l_q$, denoted by $A_1, A_2, \ldots, A_q$. 
The set of affine segments ${\cal S} = \{l_i \cap A_j: i, j = 1,2, \ldots, q\}$ is 
a parallel class of lines of $\Lambda$.
\item[(22)] The incidence structure of the set ${\cal S}$ of affine segments and 
the set of planes containing $q$ affine segments of ${\cal S}$ is isomorphic to 
an affine plane of order $q$, denoted by $\Theta({\cal S})$.
\end{enumerate}

In any projective plane of order $q$, it is well known (see, for example, \cite{Hir98})  
that the maximum value $\max{k}$ for a $k$-arc to exist is 
$$ \max{k} = \left\{
   \begin{array}{ll}
    q+1  &  \ \ \ \ \ \mbox{for $q$ odd, }\\
    q+2  &  \ \ \ \ \ \mbox{for $q$ even,}
    \end{array}
    \right.  $$
and there do exist a $(q+1)$-arc for $q$ odd and a $(q+2)$-arc for $q$ even.  

Consider the affine three space $\Lambda$ and 
the affine plane $\Theta({\cal S})$ mentioned above.
Since $\Theta({\cal S})$ is an affine plane of order $q$, 
we know that there exists in $\Theta({\cal S})$ a $(q+1)$-arc for $q$ odd 
or a $(q+2)$-arc for $q$ even.   

\begin{Theorem}
\label{exp}
If a set ${\cal K}$ of $k$ affine segments in $\Theta({\cal S})$ is a $k$-arc,  
then the set of points contained in the affine segments of ${\cal K}$ 
is a $kq$-parc in the plane $\Pi$.
\end{Theorem}

\Proof
We need only to prove that any three points not contained in a single $l_i$, 
$1 \le i \le q$,  are independent. 
Let $S_{i_1}$ and $S_{i_2}$ be two affine segments of ${\cal S}$ 
contained in $l_{i_1}$ and $l_{i_2}$, respectively, where $1 \le i_1, i_2 \le q$, $i_1 \neq i_2$.
Then there exists a unique affine Bare subplane which contains $S_{i_1}$ and $S_{i_2}$, 
and also contains one affine segment $S_i$ from each of the remaining $l_i$, $1 \le i \le q$.
This implies that any affine segment of $\cal L$ intersecting $S_{i_1}$ and $S_{i_2}$ 
intersects also each of the remaining $S_i$ for $1 \le i \le q$.
The affine segments $S_1, S_1, \ldots, S_q$ form a line of $\Theta({\cal S})$. 
Since ${\cal K}$ is a $k$-arc in $\Theta({\cal S})$, 
we know that any line of ${\cal L}$ can intersect ${\cal K}$ in at most two points, 
which means that any three points not contained in a single $l_i$, 
$1 \le i \le q$,  are independent. 
Therefore, the set of points contained in the affine segments of ${\cal K}$ is 
a $kq$-parc in the plane $\Pi$.
\qed

Let ${\cal K}'$ be a $k$-arc in a projective plane ${\Pi}_2(q)$. 
A line $l$ of $\Pi_2(q)$ is called an $i$-secant with respect to ${\cal K}'$ 
if the line $l$ intersects ${\cal K}'$ in $i$ points.  
$2$-Secants, $1$-secants and $0$-secants are sometimes
called bisecants, unisecants and external lines, respectively.
It is well known (see, for example, \cite{Hir98}) that 
there exists at least one external line for any $k$-arc.  
Let $l'$ be an external line with respect to ${\cal K}'$. 
Then ${\cal K}'$ becomes a $k$-arc in the affine plane $\Pi_2(q) \setminus l'$. 
Clearly, there exists a $(q+1)$-arc in an affine plane of odd order $q$ 
and there exists a $(q+2)$-arc in an affine plane of even order $q$.  

Let  ${\cal U} =\{U_1, U_2, \ldots, U_q\}$ be a parallel class of affine segments 
in the affine plane $\Theta({\cal S})$.
From the results in Chapter $8$ of \cite{Hir98}, 
it can be easily seen that the following structures exist 
in the affine plane $\Theta({\cal S})$ of odd order $q$:
\begin{enumerate}
\item[(23)] A $(q+1)$-arc where two affine segments of ${\cal U}$ are $1$-secants 
and $(q-1)/2$ affine segments of ${\cal U}$ are $2$-secants; 
\item[(24)] A $(q+1)$-arc where $(q+1)/2$ affine segments of ${\cal U}$ are $2$-secants; 
\item[(25)] A $q$-arc where $q$ affine segments of ${\cal U}$ are all $1$-secants.
\end{enumerate}

\begin{Theorem}
For any odd prime power $q$, there exist the following pencil arcs in $\Pi_2(q^2)$:
\begin{enumerate}
\item[(26)]  A $(q^2+q)$-parc ${\cal K}$ in which there are two lines of the pencil 
containing $q$ points of ${\cal K}$ and $(q-1)/2$ lines of the pencil containing 
$2q$ points of ${\cal K}$;
\item[(27)] A regular $((q+1)/2, 2q)$-parc;
\item[(28)] A regular $(q,q)$-parc.
\end{enumerate}
\end{Theorem}

The proof of this theorem is obvious from Theorem \ref{exp}.

The $q$ mutually disjoint affine Baer subplanes $A_1, A_2, \ldots, A_q$ 
in $\Pi \setminus l_\infty$ mentioned in (21) are incident 
with a common affine segment $G$ in $l_\infty \setminus \{P\}$. 
Also, there are other $q-1$ classes each of which contains 
$q$ mutually disjoint Baer subplanes in $\Pi \setminus l_\infty$ 
which are incident with a common affine segment in $l_\infty \setminus \{P\}$.   
These affine segments in $l_\infty$, say $G_1, G_2, \ldots, G_q$, 
are mutually disjoint in $l_\infty \setminus \{P\}$.
The affine plane $\Theta({\cal S})$ together with the set $\{ P, G_1, G_2, \ldots, G_q\}$ 
as the infinity line,  with the incidence relation of $\Pi$,  
is isomorphic to a projective plane of order $q$, denoted by $\Theta^\star({\cal S})$.  

In a projective plane of order $q$ with $q$ even, 
the following results on $k$-arcs are known:
\begin{enumerate}
\item[(29)] A $(q+1)$-arc ${\cal K}$ where $q+1$ unisecants are concurrent 
with a point called nucleus.
\item[(30)] A $(q+2)$-arc ${\cal K}$ where $(q+2)/2$ bisecants are concurrent 
in a point outside ${\cal K}$. 
\end{enumerate}

\begin{Theorem}
For any even prime power $q$,  there exist the following pencil arcs in $\Pi_2(q^2)$: 
\begin{enumerate}
\item[(31)] A regular $(q+1, q)$-parc;
\item[(32)] A regular $((q+2)/2, 2q)$-parc.
\end{enumerate}
\end{Theorem}

\Proof
We first consider a $(q+1)$-arc ${\cal K}$ in $\Theta^\star({\cal S})$ 
where the nucleus is the center $P$ of a pencil $\{l_1, l_2, \ldots, l_q, l_{\infty}\}$. 
Every line of the pencil meets ${\cal K}$ exactly once.
This implies that one affine segment in each line $l_i$, $i=1,2, \ldots, q$, 
and one affine segment in $l_\infty$ are selected. 
Then similarly to Theorem \ref{exp}, we obtain a regular $(q+1,q)$-parc in $\Pi_2(q^2)$. 

Since there exists a $(q+2)$-arc ${\cal K}'$ in $\Theta^\star({\cal S})$ 
where $(q+2)/2$ bisecants are concurrent in a point outside ${\cal K}'$,  
in a similar fashion, we can obtain a regular $((q+2)/2, 2q)$-parc in $\Pi_2(q^2)$. 
\qed

\begin{Theorem}
There exists a hierarchical arc $(q+1, q^2-q)$-harc in $\Pi_2(q^2)$ when $q$ is odd and 
a maximal $(q+2, q^2-q)$-harc when $q$ is even.
\end{Theorem}

\Proof
Let $l_\infty$ and $\Pi'$ be a line  and a Baer subplane of $\Pi_2(q^2)$, respectively,
where they intersect in $q+1$ points.
Let ${\cal K} =  {\cal K}_1 \cup  {\cal K}_2$  
be a hierarchical arc with $k_1 + k_2$ points.
${\mathcal K}_1$ is $k_1$-arc in the Baer subplane such that 
$l_\infty \cap \Pi'$ is a external line to ${\mathcal K}_1$.
Then  any  bisecant of ${\mathcal K}_1$ meet $l_\infty$ at a point
of $\Pi'$. 
Therefore we can take the points of $l_\infty \setminus \Pi'$ as
${\mathcal K}_2$.
$k_1=q+1$ when q is odd and $q+2$ when q is even are maximal points.
When q is even, $k_1 + k_2 = (q+2) + (q^2 -q) = q^2 +2$ is the maximal
value for  a hierarchical arc to exist.
\qed

\section{Conclusions}

Characterizing ideal secret sharing schemes is important in secret sharing. 
Finding an ideal secret sharing scheme for a given access structure is difficult in general. 
In this paper, we derived several combinatorial properties and 
restate the combinatorial characterization of an ideal secret sharing scheme 
in Brickell-Stinson model in terms of orthogonality of its representative array.  
We proposed practical parallel and hierarchical access structures. 
By the restated combinatorial characterization,  
we discussed sufficient conditions on finite geometries 
for ideal secret sharing schemes to realize these access structure models, 
and several infinite series of ideal secret sharing schemes 
realizing parallel or hierarchical access structure model of strength $3$ 
were also constructed from finite projective planes. 
However we understand it is a very difficult and challenging problem to 
construct ideal secret sharing schemes 
realizing parallel or hierarchical access structure model of large strengths.

\noindent
{\bf Acknowledgments:}
Research supported by Grant-in-Aid for Scientific Research (C) under Grant No.~18540109.


\end{document}